\scrollmode
\documentclass[10pt]{amsart}
\usepackage{verbatim}
\usepackage{eucal}
\usepackage{amssymb}
\usepackage{mathrsfs}
\usepackage{graphicx}
\usepackage{psfrag}
\newif \ifwide
\newif \ifavnermargin
\def \makemargins{
\ifwide
        \oddsidemargin .25in
        \evensidemargin .25in
        \textwidth 6.00in
\else
\fi
\ifavnermargin
        \headheight=7pt
        \textheight=574pt
        \textwidth=432pt
        \topmargin=14pt
        \oddsidemargin=18pt
        \evensidemargin=18pt
\else   
\fi
}

\widetrue
\makemargins

\theoremstyle{plain}
\newtheorem{theorem}[subsection]{Theorem}

\theoremstyle{definition}

\theoremstyle{remark}
\newtheorem{remark}[subsection]{Remark}

\newcount\timehh\newcount\timemm
\timehh=\time 
\divide\timehh by 60 \timemm=\time
\newcommand{\draftauthor}[1]{\author{#1
    {
      --- \protect \protect\sc\today\ ---
      \ifnum\timehh<10 0\fi\number\timehh\,:\,\ifnum\timemm<10 0\fi\number\timemm
      \protect \, \, \protect \bf DRAFT
    }
  }
}
\newcommand{\R}{{\mathbb R}}
\newcommand{\C}{{\mathbb C}}
\newcommand{\Z}{{\mathbb Z}}

\renewcommand{\P}{{\mathbb P}}

\newcommand{\fH}{{\mathfrak{H}}}

\newcommand{\TTT}{{\mathscr{T}}}

\newcommand{\SSS}{{\mathscr{S}}}
\newcommand{\EEE}{{\mathscr{E}}}
\newcommand{\MMM}{{\mathscr{M}}}
\newcommand{\ac}[1]{{\rm a.c.}\Bigl(#1\Bigr)}
\newcommand{\ee}{{\mathrm{e}}}
\newcommand{\ii}{{\mathrm{i}}}

\renewcommand{\theta}{\vartheta}

\renewcommand{\mod}{\bmod}

\DeclareMathOperator{\codim}{codim}

\DeclareMathOperator{\Proj}{Proj}
\DeclareMathOperator{\SL}{SL}

\begin{document}

\title{On toric varieties and modular forms}

\newif \ifdraft
\def \makeauthor{
\ifdraft
        \draftauthor{Paul E. Gunnells}
\else
\author{Paul E. Gunnells}
\address{Department of Mathematics\\
Rutgers University\\
Newark, NJ  07102}
\email{gunnells@andromeda.rutgers.edu}
\fi
}

\thanks{Partially supported by the NSF}

\draftfalse
\makeauthor

\ifdraft
        \date{\today}
\else
        \date{June 16, 2001}
\fi

\maketitle


\section{Introduction}\label{introduction}
Let $\ell >1$ be an integer, and consider the congruence subgroup $\Gamma _{1} (\ell ) \subset \SL_{2} (\Z)$ defined by 
\[
\Bigl (\begin{array}{cc}
a&b\\
c&d
\end{array}\Bigr)=
\Bigl (\begin{array}{cc}
1&*\\
0&1
\end{array}\Bigr)\mod \ell .
\]
Let $\MMM _{*} (\ell) = \MMM _{*} (\Gamma _{1} (\ell ), \C )$ be the
ring of holomorphic modular forms on $\Gamma _{1} (\ell )$.  In this
talk we use the combinatorics of complete toric varieties to construct
a subring $\TTT _{*} (\ell )\subset \MMM _{*} (\ell )$, the subring of
\emph{toric modular forms} (\S\ref{construction}).  This is a natural
subring, in the sense that it behaves nicely with respect to natural
operations on $\MMM _{*} (\ell )$ (namely, Hecke operators, Fricke
involution, and the theory of oldforms and newforms).  Moreover, an
explicit structure theorem for $\TTT _{*} (\ell )$ together with the
theory of Manin symbols allows us to describe the cuspidal part of
$\TTT _{*} (\ell )$ in terms of nonvanishing of special values of
$L$-functions (\S\ref{values}).  Finally, we discuss an explicit
scheme-theoretic embedding of the modular curve $X_{1} (\ell ) =
\Gamma _{1} (\ell ) \backslash \fH ^{*}$ in a weighted projective
space that was inspired by the structure of $\TTT _{*} (\ell )$
(\S\ref{equations}).

The results of \S\S\ref{construction}--\ref{values} are joint work
with Lev Borisov, and can be found in the papers
\cite{vanish,higher,toric}; the embedding of the modular curve in
\S\ref{equations} is joint work with with Lev Borisov and Sorin
Popescu and appears in \cite{modular}.  It is a pleasure to thank them
for their stimulating and interesting collaboration.


\section{Construction and basic properties}\label{construction}

Let $d$ be a positive integer, let $N$ be a rank $d$ lattice, let $M$
be the dual lattice, and let $\langle\phantom{a},\phantom{a}\rangle
\colon M\times N \rightarrow \Z $ be the pairing.  Let $\Sigma \subset
N_{\R} = N\otimes \R$ be a complete rational polyhedral fan.  A
\emph{degree function} $\deg \colon N \rightarrow \C $ is a
piecewise-linear function that is linear on the cones of $\Sigma $.
We define a function $f_{N,\deg} \colon \fH \rightarrow \C $ by
\begin{equation}\label{eff}
f_{N,\deg}(q) := \sum_{m\in M}\Bigl(\sum_{C\in\Sigma}(-1)^{\codim C}
\ac{\sum_{n\in C} q^{\langle m,n\rangle} \ee^{2\pi\ii \deg(n)}}\Bigr).
\end{equation}
Here $q=\ee ^{2\pi \ii \tau }$, where $\tau \in \fH$, the upper
halfplane, and the a.c. denotes analytic continuation.  

\begin{theorem}
Suppose $\deg$ takes values in $\ell ^{-1} \Z $, and that $\deg $ is
not integral valued on the primitive generator of any $1$-cone of
$\Sigma $.  Then $f_{N,\deg} (q)\in \MMM _{d} (\ell )$, i.e. $f$ is
the $q$-expansion of a holomorphic modular form on $\Gamma _{1} (\ell )$. 
\end{theorem}

To prove this theorem, one begins by showing that this series is
well-defined, in the sense that only finitely many terms contribute to
a given power of $q$.  To prove modularity, let $X_{\Sigma }$ be the
toric variety associated to $\Sigma $.  Then if $X_{\Sigma }$ is
smooth, one uses the Hirzebruch-Riemann-Roch theorem to show
\[
f_{N,\deg}(q)=\int_{X_{\Sigma }} \prod_D
\frac{(D/2\pi\ii)\theta(D/2\pi\ii-\alpha _{D},\tau )\theta'(0,\tau
)}{\theta(D/2\pi\ii,\tau )\theta(-\alpha_{D},\tau )},
\]
where $D$ ranges over the torus-invariant divisors of $X_{\Sigma }$,
$\alpha _{D}$ is the value of $\deg$ on the primitive generator of the
$1$-cone corresponding to $D$, and $\theta (z,\tau )$ is Jacobi's
theta function.  The expression in the integrand is evaluated in the
cohomology ring $H^{*} (X_{\Sigma })$ using the triple product formula
for $\theta $.   The case of singular $X_{\Sigma }$ is handled using a
limiting argument.
       
Let $\TTT _{*} (\ell )$ be the full subring of $\MMM _{*} (\ell )$
given by taking all $\C $-linear combinations of all $f_{N,\deg}$.  
\begin{theorem}
$\TTT _{*} (\ell )$ is closed under the action of the Hecke operators,
the Fricke involution, and Atkin-Lehner lifting.  
\end{theorem}
These statements are proved by direct manipulations of fans and degree
functions.  For example, if $p$ is a prime not dividing $\ell $, the
action of Hecke operator $T_{p}$ on $f_{N,\deg}$ is given by
\[
T_{p} f_{N,\deg} = \frac{p-p^{d-1}}{p-1}f_{N,\deg} + \sum _{S} f_{S,p\deg},
\]
where the sum is taken over lattices $S$ satisfying $N\subset S\subset
\frac{1}{p}N$ and $[S:N] = p^{d-1}$.  Note that these lattices are
\emph{not} the usual lattices involved in the definition of the Hecke
operators.

\begin{remark}
The definition \eqref{eff} was motivated by L. Borisov and A. Libgober's
computation of the elliptic genera of toric varieties \cite{borlibg}.
Also, similar sums were studied by W. Nahm, who showed that
they had (quasi)modular properties \cite{nahm}.
\end{remark}


\section{Special values of $L$-functions}\label{values}
A natural question is how far is $\TTT _{*} (\ell )$ from being all of
$\MMM _{*} (\ell )$.  The inclusion $\TTT_{*} (\ell ) \subset \MMM_{*}
(\ell )$ is certainly proper, 
since for example there are no weight 1 cusp forms in
$\TTT _{1} (\ell )$.   However, it turns out that almost all
modular forms are toric.  To state the precise result, let
$\xi=\ee^{2\pi\ii/\ell }$, and 
for $0<a<\ell $
let $s_{a} = s_{a} (q)$ be the weight 1 Eisenstein series 
\[
s_{a} (q) = \frac1{2\pi\ii}\frac d{dz} \ln \theta(z,\tau)|_{z=a/\ell
}= \frac{\xi^a+1}{2(\xi^{a}-1)} - \sum _{d}q^{d}\sum _{k|d} (\xi^{ka}-
\xi^{-ka}).
\]
For each $k$, let $\MMM _{k} (\ell ) = \SSS_{k} (\ell )\oplus \EEE
_{k} (\ell )$ be the decomposition into cusp forms and Eisenstein
series.  Given a Hecke eigenform $f$, let $L (f,s)$ be the associated
$L$-function.

\begin{theorem}
The ring $\TTT _{*} (\ell )$ is multiplicatively generated by the
$s_{a}$, $0<a<\ell $.  In weight two, $\TTT _{2} (\ell )$ modulo $\EEE
_{2} (\ell )$ is equal to the $\C $-span of all Hecke eigenforms $f$
with $L (f,1) \not = 0$.  For weights $ k\geq 3$, the space $\TTT _{k}
(\ell )$ coincides with $\MMM _{k} (\ell )$ modulo $\EEE _{k} (\ell
)$.
\end{theorem}

In other words, the cuspidal part of $\TTT _{k} (\ell )$ is easy to
describe: for weights $\geq 3$ all cusp forms are toric, and for
weight $2$ only those in the span of the ``analytic rank 0'' forms are
toric.  In general, however, it is not clear what Eisenstein series
are toric.  For example, in weight $2$ the space $\TTT _{2} (25) \cap
\EEE _{2} (25)$ has codimension one in $\EEE _{2} (25)$.

The proof of this theorem is a computation with \emph{Manin symbols}
\cite{merel}.  The space $M (\ell )$ of Manin symbols of
level $\ell $ is the $\C $-vector space generated by the symbols 
\[
\left\{(a,b) \in (\Z /\ell \Z )^{2}\mid \Z a + \Z b = \Z /\ell \Z\right\},
\]
modulo certain relations.  This space is dual to the space $\SSS _{2}
(\ell )$ in the following sense.  To each symbol $(a,b)$, one can
associate an ideal geodesic $\gamma $ on $\fH$.  There is an
involution on $M (\ell )$ that defines two subspaces $M^{\pm } (\ell
)$.  Then there are two subspaces $S^{\pm } (\ell ) \subset M^{\pm }
(\ell )$ such that the pairings $S^{\pm }(\ell ) \times \SSS _{2}
(\ell ) \rightarrow \C $ given by integration of a cuspform along a
chain are perfect.  Moreover the two subspaces $S^{\pm } (\ell )$ are
dual to each other via the intersection pairing on cycles on $X_{1}
(\ell )$.

The key point to the computation is that modulo Eisenstein series the
product $s_{a}s_{b}$ has properties similar to the the symbol
$(a,b)\in M^{-} (\ell )$.  This allows us to define a map $\mu \colon
M^{-} \rightarrow \MMM _{2} (\ell )/\EEE _{2} (\ell ) \cong \SSS _{2}
(\ell ) $ by $(a,b) \mapsto s_{a}s_{b}$.  Then we define a map
$\varphi \colon \SSS _{2} (\ell )\rightarrow \SSS _{2} (\ell )$ by
\[
f = \sum a_{n}q^{n} \longmapsto \sum a_{n} \Bigl(\int _{0}^{\ii \infty }
(T_{n}f)\,d\tau\Bigr)  q^{n},
\]
where $T_{n}$ is the $n$th Hecke operator.
This map has the property that $\varphi (f) = 0$ if $f$ is a new
eigenform with $L (f,1)=0$.  Finally using Merel's universal formula for the
Hecke action on Manin symbols, an explicit description of the
intersection pairing on $X_{1} (\ell )$, and some manipulations with
the product of $q$-expansions $s_{a} (q) s_{b} (q)$, we show that
$\varphi $ factors through $\mu $.  

For higher weights the argument is similar and its conclusion is
identical.  We then appeal to a theorem of Jacquet and Shalika
\cite{js} that implies $L (f,1)$ cannot vanish for a new eigenform of
weight $\geq 3$.


\section{Equations of modular curves}\label{equations}
Let $p\geq 5$ be a prime.
Given the simple description of $\TTT _{*} (p)$ and its relation
with the forms of analytic rank 0, one is interested in estimating $\dim
\TTT _{2} (p)$.  As a first step in this direction, we have studied
the (incomplete) linear system on the modular curve $X_{1} (p)$
induced the weight 1 Eisenstein series $s_{a}$.

\begin{theorem}
Let $\P $ be the weighted projective space 
\[
\Proj \C [s_{a}, t_{b} \mid 0 < a,b < p ],
\]
where $\deg s_{a} = 1$, and $\deg t_{b} = 2$.  Then the modular curve
$X_{1} (p)$ is scheme-theoretically cut out from $\P $ by the equations (1)
$s_{a} = - s_{p-a}$, (2) $t_{b} = t_{p-b}$, and (3) $s_{a}s_{b} +
s_{b}s_{c} + s_{c}s_{a} = t_{a} +t_{b} +t_{c}$ if $a+b+c = 0\mod p$.
\end{theorem}
Here the variables $t_{b}$ correspond to certain weight 2 Eisenstein
series on $X_{1} (p)$.

To prove this theorem, we first show that the map $X_{1} (p)
\rightarrow \P ^{(p-3)/2}$ defined by $\tau \mapsto \{s_{a} (\tau )
\}$ is a closed embedding.
Then we construct a system of
differential equations 
\[
\frac{dr_{a}}{dz} = -\frac{1}{p-2}\sum _{k\not =0,a} r_{k}r_{a-k} +
2r_{a}s_{a}, \quad a \in (\Z /p\Z )^{\times }  
\]
that mimics the system satisfied by the set of elliptic functions 
\[
z \longmapsto \frac{\theta (a/p-z,\tau )\theta _{z} (0,\tau )}{\theta
(-z,\tau )\theta (a/p, \tau )}, 
\]
having poles only along a $p$-torsion subgroup.  We construct
``standard solutions''
\[
r_{a} (z) = \frac{1}{z} +s_{a} +t_{a}z + \cdots,
\]
which satisfy certain quadratic relations.  If we define a ring using
these relations, we get a ring with the same Hilbert function as the
coordinate ring of an elliptic normal curve $C$ of degree $p$ in $\P
^{p-1}$.  Then we show that $C$ is singular if and only if $C$ is a
$p$-gon if and only if the coordinates $\{s_{a},t_{b} \}$ correspond
to a cusp of $X_{1} (p)$.  Finally we show that a deformation of a
solution to our system leads to a deformation of the elliptic curve,
which leads to an identification of the scheme defined by the
equations in the theorem with $X_{1} (p)$.



\providecommand{\bysame}{\leavevmode\hbox to3em{\hrulefill}\thinspace}

\end{document}